\documentclass[12pt, letterpaper, twoside]{article}
\usepackage[utf8]{inputenc}
\usepackage{amsfonts,amssymb,amsthm,amsmath}

\theoremstyle{plain}
\newtheorem{theorem}{Theorem}

\theoremstyle{definition}

\theoremstyle{remark}
\newtheorem*{remark}{Remark}

\newcommand{\N}{\mathbb{N}}
\newcommand{\R}{\mathbb{R}}

\newcommand{\de}{\mathrm{d}}

\newcommand{\id}{\operatorname{id}}

\newcommand{\fuermich}[1]{}
\newcommand{\ifindoubt}[1]{}
\newcommand{\suppressed}[1]{}

\fuermich{}

\usepackage{enumitem}

\title{Nowhere vanishing primitive of a symplectic form}
\author{B. Stratmann}

\date{}

\begin{document}
\maketitle
\begin{abstract}
Let $M$ be a manifold with an exact symplectic form~$\omega$. Then there is a nowhere
vanishing primitive $\beta$ for $\omega$, i.e. $\omega=\de\beta$.
\end{abstract}
Blohmann and Weinstein 
investigate the tangent Lie algebroid $TM$ 
over a symplectic manifold $(M,\omega)$ and they
are able to prove that it admits a Hamiltonian structure if and only if the symplectic form $\omega$ is exact
and there is a nowhere vanishing primitive. In the sequel they
raise the question whether each
manifold with an
exact symplectic form admits a nowhere vanishing primitive, i.e. if there
exists a nowhere vanishing 1-form $\beta$
with $\omega=\de\beta$
(\cite[Cor. 6.14]{BloWei18}). Here we give
an affirmative answer.
\begin{theorem}
Let $M$ be a manifold with an exact symplectic form $\omega$. Then there is a nowhere
vanishing primitive $\beta$ for $\omega$.
\label{theorem_symplectic}
\end{theorem}
\begin{remark}
As an immediate consequence a manifold with an exact symplectic form admits a nowhere vanishing Liouville vector field.
\end{remark}
\begin{proof}
Without restriction we may assume $M$ to be connected. Denote by $2n$ the dimension of $M$ and fix an orientation.
Exactness of $\omega$ implies that $M$ is non-compact.\\
Choose a Morse function $\gamma:M\to\R^{>0}$
which is exhaustive, i.e. the sets
$\{x\in M\ \vert\ \gamma(x)\leq c\}$ are compact for all $c\in\R$.
Denote $C_0$ the discrete set of critical points of $\gamma$.\\
The first step is to construct a primitive $\theta_1$ of $\omega$ which is of a particular form
around each point $z\in C_0$.
For this choose for each
$z\in C_0$ a contractible open neighborhood
$U_z$ in a way that
$U_z$ and $U_{\tilde z}$ do not intersect
for all $z\not =\tilde z$ and such that
there are $2n$ global functions $x_i$ with $x_i(z)=0$
for all $z\in C_0$ providing local coordinates
on $U_z$ such that
\begin{equation}
 \omega=2\cdot\sum_{i=1}^n
 \de x_i\wedge\de x_{i+n}\qquad\text{on}\ U_z\ .
\end{equation}
By assumption there is a primitive $\theta_0$
for $\omega$.
Locally on each neighborhood $U_z$ the 1-form
$$\tau_z=\sum_{i=1}^n x_i\de x_{i+n}-\sum_{i=1}^n
 x_{i+n}\de x_i$$
is a primitive as well.
As $U_z$ is contractible there is a function
$f_z:U_z\to \R$ such that
\begin{equation}
 \theta_0-\tau_z=\de f_z\qquad\text{on}\ U_z\ .
\end{equation}
Choosing smooth functions $\chi_z:M\to [0,1]$
with compact support in $U_z$ and a
relatively compact open neighborhood $V_z\subset U_z$ of $z$ on which $\chi_z$ is identically equal 
to $1$
we define the 1-form
\begin{equation}
\theta_1=\theta_0-\sum_{z\in C}\de(\chi_z\cdot f_z)\ ,
\end{equation}
again a primitive of $\omega$. By construction 
$\theta_1$ and $\tau_z$ are equal on $V_z$.\\
The next step will be to modify $\gamma$ to $\nu$
such that the function $\nu$
has the same critical points but is of a particular form in the neighborhood of each critical point with respect to the 
coordinate functions $x_i$ already introduced.
This will be done by defining a diffeomorphism $\Psi:M\to M$ fixing $C_0$ and
setting $\nu=\gamma\circ\Psi$.\\
As $\gamma$ is a Morse function, the Morse Lemma states that, after possible shrinking of $V_z$,
there are $2n$
global functions $y_i$
with $y_i(z)=0$ defining local 
coordinates on $V_z$ for all $z\in C_0$
such that there are integers $0\leq s_z\leq 2n$ so that
\begin{equation}
\gamma(y)=\sum_{i=1}^{s_z}y_i^2-\sum_{i=s_z+1}^{2n}y_i^2+\gamma(z)\qquad\text{on}\ V_z
\end{equation}
This choice can be done in a way that the two bases
$\de x_1,\ldots,\de x_{2n}$ and 
$\de y_1,\ldots,\de y_{2n}$ in $T^*_zM$ induce the same chosen orientation, e.g. just by replacing $y_1$ by $-y_1$ locally in $V_z$ if necessary.
After possible further shrinking of $V_z$
there is an orientation preserving diffeomorphism
$\psi_z:V_z\to\psi_z(V_z)$
with $\psi_z^*y_i=x_i$ for all $1\leq i\leq 2n$
and $\psi_z(z)=z$. 
It is a classical result that there is a
relatively compact open neighborhood $V_z'\subset V_z$ of $z$
and an orientation preserving diffeomorphism
$\Psi_z:M\to M$ such that
$\Psi_z\vert_{V_z'}=\psi_z\vert_{V_z'}$
and $\Psi_z\vert_{M\backslash V_z}
=\id\vert_{M\backslash V_z}$. The result may be deduced from \cite[Theorem~5.5]{Pal59} when applied to $V_z'$.
So there is a global diffeomorphism $\Psi:M\to M$
defined by
\begin{equation*}
 \Psi(x)=\begin{cases}
          \Psi_z(x)&\text{for}\quad x\in V_z\\
          x&\text{elsewhere}
         \end{cases}
\end{equation*}
by which we define the exhaustive
Morse function
$\nu=\gamma\circ\Psi$ which has the same set $C_0$ as set of critical points.\\
In summary, we constructed an exhaustive
Morse function $\nu$ with
discrete set $C_0$ of
critical points and not intersecting relatively
compact open neighborhoods $V_z'$
such that for each $z\in C_0$ 
\begin{equation}
\de\nu=2\cdot\sum_{i=1}^{s_z}x_i\de x_i-2\cdot\sum_{i=s_z+1}^{2n}x_i\de x_i\qquad\text{on}\ {V_z'}
\label{equation_localFormRho}
\end{equation}
and a primitive $\theta_1$ which satisfies
\begin{equation}
\theta_1=\sum_{i=1}^n x_i\de x_{i+n}-
\sum_{i=1}^n x_{i+n}\de x_i\qquad\text{on}\ V_z'\ .
\label{equation_localFormTheta}
\end{equation}
There is a $K_0>0$ such that for all $K\geq K_0$
the system of linear equations
\begin{equation}
\begin{split}
\left(\sum_{i=1}^n x_i\de x_{i+n}-
\sum_{i=1}^n x_{i+n}\de x_i\right)
\qquad\qquad\mbox{ }\qquad\\
+\  K\cdot\left(2\cdot\sum_{i=1}^{s}x_i\de x_i-2\cdot \sum_{i=s+1}^{2n}x_i\de x_i\right)=0
\end{split}
\end{equation}
has $(x_1,\ldots,x_{2n})=(0,\ldots,0)$
as the only solution for all $1\leq s\leq 2n$.
Applying this to the equations~\eqref{equation_localFormRho}
and \eqref{equation_localFormTheta}
we obtain for all $z\in C_0$
\begin{equation*}
\left\{
x\in{V_z'}\ \vert\ \theta_1(x)+K\cdot\de\nu(x)=0\right\}=\{z\}\qquad\text{for all}\ K\geq K_0
\end{equation*}
and as $K_0$ is independent of $z\in C_0$
\begin{equation*}
\bigcup_{z\in C}\left\{
x\in{V_z'}\ \vert\ \theta_1(x)+K\cdot\de\nu(x)=0\right\}=C\qquad\text{for all}\ K\geq K_0\ .
\end{equation*}
So we constructed an exhaustive 
function $\nu$ such that for all $K\geq K_0$
the primitive
$$\theta_1+ K\cdot \de \nu$$
vanishes inside $V'=\bigcup_{z\in C}V_z'$
exactly on $C_0$.\\
The next step will be 
to define a primitive which does not vanish outside $V'$ at all and inside $V'$ still exactly on $C_0$. For this 
we modify $\nu$ by a smooth function $\lambda:\R^{>0}\to \R^{>0}$
defining $\rho=\lambda\circ\nu$
such that the primitive
\begin{equation}
 \beta_0=\theta_1 + \de\rho
\end{equation}
vanishes exactly at $C_0$.\\
Recall that on the closed subset 
$M\backslash V'$ the 1-form $\de\nu$ vanishes nowhere.
For each $m=1,2,\ldots$
the set
$W_m=\nu^{-1}([m-1,m])\backslash V'$
is compact using the fact that $\nu$
is exhaustive.
So there is a $K_m>K_0$ such that
\begin{equation*}
\{w\in W_m\ \vert \ \theta_1(w)+K\cdot\de\nu(w)=0\}
\end{equation*}
is empty for all $K>K_m$.
There is a smooth function
$\lambda:\R^{>0}\to\R^{>0}$
with $\lambda'(s)>K_m$ for all $s\in [m-1,m]$.
Now with $\rho=\lambda\circ\nu$ we get
\begin{equation}
 \de\rho(x)=\lambda'(\nu(x))\cdot\de\nu(x)
\end{equation}
with $\lambda'(\nu(x))$ chosen sufficiently large that
$\beta_0=\theta_1 +\de\rho$ 
vanishes on $V'$ exactly on $C_0$ while
outside of $V'$ it does
not vanish at all.\\
So there is a primitive $\beta_0$ with discrete vanishing set $C_0$. 
The final step now consists in moving away these zeroes successively by induction. 
For this we will construct a sequence of compact subsets $Z_m\subset M$ 
with the following properties. The sequence
is exhaustive, i.e. $\bigcup_{m\in\N}Z_m=M$, and
for each $m\in\N$ the compact set
$Z_m$ is contained in the interior of
$Z_{m+1}$ and 
each point in $Z_{m+1}\backslash Z_m$ is joint by a smooth path in $M\backslash Z_m$ to a point in $M\backslash Z_{m+1}$.
Further the sequence is initialized such that $Z_0\cap C_0=\emptyset$. In order to
define $Z_m$
choose a smooth exhaustion function
$\kappa:M\to\R^{>0}$ and a strictly increasing sequence $v_m\in\R^{>0}$ of regular values of $\kappa$
such that $v_m$ converges to infinity and
$\{x\in M\ \vert\ \kappa(x)\leq v_0\}$
is non-empty, disjoint to $C_0$ and has connected complement.
Set $Z_m'=\{x\in M\ \vert\ \kappa(x)\leq v_m\}$
and let $D_{m 1},\ldots, D_{m j_m}$
be those connected components of $M\backslash Z_m'$ whose closure in $M$ is compact. Note that compactness of $Z_m'$ implies
that there are only finitely many sets $D_ {m j}$ for fixed $m$. Hence the sets
$Z_m=Z_m'\cup\bigcup_{j=1}^{j_m}\overline{D_{m j}}$ define the exhausting sequence with the desired properties.\\
We will construct an
adapted sequence of 1-forms $\beta_m$ with vanishing sets
$C_m=\{x\in M\ \vert \ \beta_m(x)=0\}$
such that
\begin{enumerate}[label=(\arabic*)]
\item $\omega=\de\beta_m$
\item $C_m$ is discrete
\item $Z_m\cap C_m=\emptyset\quad$ and
\item $\beta_{m}\vert_{Z_{m-1}}=\beta_{m-1}\vert_{Z_{m-1}}\quad$ for all $m\geq 1$\label{item_betaLocallyStable}
\end{enumerate}
holds. We may consider Condition~\ref{item_betaLocallyStable} to be empty for $m=0$ such that $\beta_0$ satisfies the four conditions. For the induction step
consider a 1-form $\beta_m$ to
be given for some $m \in\N$ satisfying the four conditions.
Since $Z_{m+1}$ is compact, the set
$D_{m+1}=Z_{m+1}\cap C_m$ is finite.
By assumption
each point in $Z_{m+1}\backslash Z_m$ is joint by a smooth path in
$M\backslash Z_m$ to a point in $M\backslash Z_{m+1}$.
So we may choose a
relatively compact open subset
$U_m\subset M\backslash Z_m$ with
$U_m\cap C_m=D_{m+1}$ such that the intersection of each connected component of $U_m$ with $(M\backslash Z_{m+1})$
is not empty.
There is a symplectomorphism of $U_m$
which maps $D_{m+1}$ to a set contained in $M\backslash Z_{m+1}$ and this symplectomorphism
can be chosen such that it differs from the identity only on a relatively compact subset of $U_m$ (\cite{Boo69}). Thus the symplectomorphism can be extended as the identity on $M\backslash U_{m}$ to a symplectomorphism
$\varphi_m:M\to M$. We may define
$\beta_{m+1}=\varphi_m^*\beta_m$.
By construction $C_{m+1}=\varphi_m^{-1}(C_m)\subset M\backslash Z_{m+1}$ since $C_{m+1}\backslash U_m=C_m\backslash U_m$
and $\varphi_m\vert_{M\backslash U_m}=\id\vert_{M\backslash U_m}$.\\
As the sequence $\beta_m$ becomes stationary
on some neighborhood of each point it converges to
a 1-form $\beta$ with $\omega=\de\beta$.
Since $\beta$ does not vanish in $Z_m$
for all $m$, the form $\beta$ is a nowhere
vanishing primitive for $\omega$.
\end{proof}
\paragraph{Acknowledgement}
We thank P.~Heinzner for his ongoing support
over the last years. Further we have to thank A.~Abbondandolo and
M.~Kalus for helpful remarks and comments
and L.~Ryvkin for directing our attention to the question. 

\bibliographystyle{alpha}


\end{document}